\documentstyle[amscd,amssymb,verbatim,epsf,epsfig,latexsym,float,subfigure]{amsart}

\begin{document}

\theoremstyle{plain}
\newtheorem{theorem}{Theorem}
\newtheorem{Cor}{Corollary}
\newtheorem{Con}{Conjecture}
\newtheorem{Main}{Main Theorem}
\newtheorem{lemma}{Lemma}
\newtheorem{proposition}{Proposition}
\newenvironment{proof}{{\bf Proof:} }{\hfill $\Box$
\mbox{}}

\theoremstyle{definition}
\newtheorem{definition}{Definition}[section]
\newtheorem{Note}{Note}

\theoremstyle{example}
\newtheorem{example}{Example}[section]

\theoremstyle{remark}
\newtheorem{remark}{Remark}[section]

\theoremstyle{remark}
\newtheorem{notation}{Notation}
\renewcommand{\thenotation}{}

\errorcontextlines=0
\numberwithin{equation}{section}
\renewcommand{\rm}{\normalshape}%
\newcommand{\ici}[1]{\stackrel{\circ}{#1}}
\title[The use of soft matrices on soft multisets in an optimal decision process]%
   {The use of soft matrices on soft multisets in an optimal decision process}

\author[Arzu Erdem, Cigdem Gunduz Aras, Ayse Sonmez, and H\"usey\.{I}n \c{C}akall\i]{Arzu Erdem*, Cigdem Gunduz Aras* , Ayse Sonmez**, and H\"usey\.{I}n \c{C}akall\i\;***\\
*Kocaeli University, Department of Mathematics, Kocaeli, Turkey Phone:(+90262)3032102\\ **Gebze Institute of Technology, Department of Mathematics,  Gebze-Kocaeli, Turkey Phone: (+90262)6051389 \\*** Maltepe University, Marmara E\u{g}\.{I}t\.{I}m K\"oy\"u, TR 34857, \.{I}stanbul-Turkey \; \; \; \; \; Phone:(+90216)6261050 ext:2248, \;  fax:(+90216)6261113 }

\address{Arzu Erdem, Kocaeli University, Department of Mathematics, Kocaeli, Turkey Phone:(+90 262) 303 2102}
\email{erdem.arzu@@gmail.com}
\address{Cigdem Gunduz Aras, Kocaeli University, Department of Mathematics, Kocaeli, Turkey Phone:(+90 262) 303 2102}
\email{carasgunduz@@gmail.com; caras@@kocaeli.edu.tr}
\address{Ay\c{s}e S\"onmez, Department of Mathematics, Gebze Institute of Technology, Cayirova Campus 41400 Gebze- Kocaeli, Turkey Phone: (+90 262) 605 1389}
\email{asonmez@@gyte.edu.tr; ayse.sonmz@@gmail.com}
\address{Huseyin \c{C}akall\i\;Maltepe University, Department of Mathematics, Marmara E\u{g}\.{I}t\.{I}m K\"oy\"u, TR 34857, Maltepe, \.{I}stanbul-Turkey \; \; \; \; \; Phone:(+90 216) 626 1050 ext:2248, \;  fax:(+90 216) 626 1113}
\email{hcakalli@@maltepe.edu.tr; hcakalli@@gmail.com}

%\author{}


\date{\today}

\maketitle

\begin{abstract}
In this paper, we introduce a concept of a soft matrix on a soft multiset, and investigate how
to use soft matrices to solve decision making problems. An algorithm for a multiple choose selection
problem is also provided. Finally, we
demonstrate an illustrative example to show the decision making steps.
\end{abstract}

\maketitle

\section{Introduction}

\normalfont{}
One has to make a choose between alternative actions in almost every days
and make decisions. Most decisions involve multiple objectives. For example, a developer and manufacturer of computer equipment,
the X Industry  would
like to manufacture a large quantity of a product if the consumers demand
for the product adequately high. Unfortunately, sometimes the technology
development is difficult and there are some crucial elements such as target
time for getting market, competitors, etc.  The Development Group of the X
industry argues that an untested product is introduced and they propose
moving up the production start,  and putting the device into production
before final tests, meanwhile launching a high-priced advertising campaign
offering the units as available now. This shows a decision event activity.
When we mention decision making one needs to figure out what to do in the
face of difficult circumstances. Molodtsov \cite{Molodtsov} has introduced
`Soft Set Theory' in which we can use parametrization tools for dealing
with this kind of uncertainty and difficulties in the process of decision
making. Maji,~Biswas~and~Roy \cite{Maji} defined equality of two soft sets,
a subset and a super set of a soft set, complement of a soft set, null soft set,
and absolute soft set with examples. In the study, soft binary operations
like AND, OR and the operations of union, intersection were also
characterized. Sezgin and Atag\"{u}n \cite{Sezgin} proved that certain De
Morgan's law hold in soft set theory with respect to different operations on
soft sets. Ali, Feng, Liu, Min and Shabir \cite{Ali} introduced some new
notions such as the restricted intersection, the restricted union, the
restricted difference and the extended intersection of two soft sets. The
relationship among soft sets, soft rough sets and topologies was established
by Li and Xie in \cite{Li}. Based on the novel granulation structures called
soft approximation spaces, soft rough approximations and soft rough sets
were introduced by Feng, Liu, Fotea and Jun in \cite{Feng}. Jiang, Tang, Chen,
Liu and Tang presented an extended fuzzy soft set theory by using the
concepts of fuzzy description logics to act as the parameters of fuzzy soft
sets. Gunduz and Bayramaov in \cite{Caras} introduced some important properties of fuzzy soft topological spaces. An alternative approach to attribute reduction in multi-valued
information system under soft set theory was presented by Herawan, Ghazali
and Deris in \cite{Herawan}. Gunduz, Sonmez and Cakalli in \cite{Caras1} introduced  soft open and soft closed mappings, soft homeomorphism concept. An application of soft sets to a decision making
problem with the help of the theory of soft sets was studied in  \cite{Cagman, Han, Maji2, Mamat, Singh, Zhang}.

When we consider a multiple objective decision making problem with $m$
criteria and $n$ alternatives, it is easy to show the decision making
methodology as a decision table. This representation has several advantages.
It is easy to store and manipulate matrices and hence the soft sets
represented by them in a computer. Soft matrices which were a
matrix representation of the soft sets were introduced in \cite{Basu,Cagman2,Cagman3,Mondal2,Mondal,Vijayabalaji}.
However, some decision making problems include several decision-makers such
as a group decision making problem. As a generalization of Molodtsov's soft
set, the definition of a soft multiset was introduced by Alkhazaleh and
Salleh in \cite{Alkhazaleh}. In the research, they gave basic operations
such as complement, union and intersection with examples and then in \cite{Alkhazaleh2}\ they introduced the definition of fuzzy soft multiset as a
combination of soft multiset and fuzzy set and study its properties and
operations.

In this paper we improve a soft multiset concept in decision-making problems
and apply soft matrices concept to group decision making problems.

\maketitle

\section{Preliminaries}
First of all, we recall some basic concepts and notions, which are necessary
foundations of group decision making methods.

\begin{definition}
(\cite{Molodtsov}) Let $U$ be an initial universe, $P\left( U\right) $ be
the power set of $U$, $E$ be a set of all parameters and $A\subset E$. A
soft set $\left( f_{A},E\right) $ on the universe $U$ is defined by the set
of ordered pairs%
\begin{equation*}
\left( f_{A},E\right) :=\left\{ \left( e,f_{A}\left( e\right) \right) :e\in
E,~f_{A}\left( e\right) \in P\left( U\right) \right\}
\end{equation*}%
where $f_{A}:E\rightarrow P\left( U\right) $ such that $f_{A}\left( e\right)
=\emptyset $ if $e\notin A$.
\end{definition}

\begin{example}
Let us consider a soft set $\left( f_{E},E\right) $ which describes the
\textquotedblright color of the shirts\textquotedblright\ that Mrs. X is
considering to buy. Suppose that there are five shirts in the universe $%
U=\{s_{1},s_{2},s_{3},s_{4},s_{5}\}$ under consideration, and that $%
E=\{e_{1},e_{2},e_{3},e_{4}\}$ is a set of decision parameters. For each $%
e_{i},$ $i=1,2,3,4,$ denotes the parameters "white", "purple", "red" and
"blue", respectively. Let $A=\{e_{1},e_{3}\}\subset E$ and
\begin{eqnarray*}
f_{A}\left( e_{1}\right) &=&U, \\
f_{A}(e_{3}) &=&\{s_{2},s_{4}\}.
\end{eqnarray*}%
Then we can view the soft set $\left( f_{A},E\right) $ as consisting of the
following collection of approximations:%
\begin{equation*}
\left( f_{A},E\right) =\left\{ \left( e_{1},U\right) ,\left(
e_{3},\{s_{2},s_{4}\}\right) \right\} .
\end{equation*}
\end{example}

\begin{definition}
(\cite{Alkhazaleh}) Let $\{U_{i}:i\in I\}$ be a collection of universes such
that $\cap _{i\in I}U_{i}=\emptyset $, $\left\{ E_{i}=E_{U_{i}}:i\in
I\right\} $ be a collection of sets of parameters, $E=\prod\limits_{i\in
I}E_{i},$ $U=\prod\limits_{i\in I}P\left( U_{i}\right) ,$ $A\subset E.$ A
pair $(F_{A},E)$ is called a soft multiset over $U$,where $F_{A}$ is a
mapping given by $F_{A}:A\rightarrow U.$
\end{definition}

This paper will focus on the situation that universe sets $U_{i}$ and
parameter sets $E_{i}$ are both finite sets for each $i\in I$ .

\begin{example}
\label{ex2}Suppose that there are three universes $U_{1},U_{2},U_{3}.$ Let
us consider a soft multiset $(F_{A},E)$ which describes the "attractiveness
of houses", "attractiveness of cars"\ and "attractiveness of hotels"\ that
Mrs. X is considering for accommodation purchase, transportation purchase,
and venue to hold a wedding celebration respectively.%
\begin{eqnarray*}
U_{1} &=&\{h_{1},h_{2},h_{3},h_{4},h_{5},h_{6}\} \\
U_{2} &=&\{c_{1},c_{2},c_{3},c_{4},c_{5}\} \\
U_{3} &=&\{v_{1},v_{2},v_{3},v_{4}\} \\
E_{1} &=&E_{U_{1}}=\{e_{11}="expensive",e_{12}="cheap",e_{13}="4\text{ }%
bedroom\text{ }flat", \\
e_{14}&=&"3\text{ }bedroom\text{ }and\text{ }terraced%
\text{ }house", e_{15} ="located\text{ }in\text{ }the\text{ }heart\text{ }of\text{ }the%
\text{ }city"\} \\
E_{2} &=&E_{U_{2}}=\{e_{21}="expensive",e_{22}="cheap",e_{23}="friendly\text{
}technology", \\
e_{24} &=&"better\text{ }performance",e_{25}="luxury",e_{26}="Made\text{ }in%
\text{ }Germany"\} \\
E_{3} &=&E_{U_{3}}=\{e_{31}="expensive",e_{32}="cheap",e_{33}="in~\dot{I}%
stanbul",
\\
e_{34}&=&"located\text{ }in\text{ }the\text{ }historic\text{ }centre", e_{35} ="neoclassic\text{ }hotel"\} \\
E &=&\prod\limits_{i\in I}E_{i},U=\prod\limits_{i\in I}P\left( U_{i}\right)
\\
A &=&\{a_{1}=\left( e_{11},e_{21},e_{31}\right) ,a_{2}=\left(
e_{11},e_{22},e_{34}\right) ,a_{3}=\left( e_{12},e_{23},e_{35}\right)
,a_{4}=\left( e_{15},e_{24},e_{32}\right) , \\
a_{5} &=&\left( e_{14},e_{23},e_{33}\right) ,a_{6}=\left(
e_{12},e_{25},e_{32}\right) ,a_{7}=\left( e_{13},e_{21},e_{31}\right)
,a_{8}=\left( e_{11},e_{26},e_{32}\right) \}\subset E
\end{eqnarray*}%
Suppose that%
\begin{eqnarray*}
F_{A}\left( a_{1}\right) &=&\left(
\{h_{3},h_{4},h_{5},h_{6}\},\{c_{1},c_{2},c_{3}\},\{v_{2},v_{3}\}\right) , \\
F_{A}\left( a_{2}\right) &=&\left(
\{h_{3},h_{4},h_{5},h_{6}\},\{c_{4},c_{5}\},\{v_{1},v_{2}\}\right) , \\
F_{A}\left( a_{3}\right) &=&\left( \{h_{1},h_{2}\},\emptyset
,\{v_{2},v_{3}\}\right) , \\
F_{A}\left( a_{4}\right) &=&\left(
U_{1},\{c_{3},c_{4}\},\{v_{1},v_{4}\}\right) , \\
F_{A}\left( a_{5}\right) &=&\left( \{h_{3},h_{4},h_{5}\},\emptyset
,\{v_{2},v_{4}\}\right) , \\
F_{A}\left( a_{6}\right) &=&\left(
\{h_{1},h_{2}\},U_{2},\{v_{1},v_{4}\}\right) , \\
F_{A}\left( a_{7}\right) &=&\left( \emptyset
,\{c_{1},c_{2},c_{3}\},\{v_{2},v_{3}\}\right) , \\
F_{A}\left( a_{8}\right) &=&\left(
\{h_{3},h_{4},h_{5},h_{6}\},\{c_{4},c_{5}\},\{v_{1},v_{4}\}\right) ,
\end{eqnarray*}%
Then we can view the soft multiset $(F_{A},E)$ as consisting of the
following collection of approximations%
\begin{eqnarray*}
(F_{A},E) &=&\{\left( \left( e_{11},e_{21},e_{31}\right) ,\left(
\{h_{3},h_{4},h_{5},h_{6}\},\{c_{1},c_{2},c_{3}\},\{v_{2},v_{3}\}\right)
\right) , \\
&&\left( \left( e_{11},e_{22},e_{34}\right) ,\left(
\{h_{3},h_{4},h_{5},h_{6}\},\{c_{4},c_{5}\},\{v_{1},v_{2}\}\right) \right) ,
\\
&&\left( \left( e_{12},e_{23},e_{35}\right) ,\left(
\{h_{1},h_{2}\},\emptyset ,\{v_{2},v_{3}\}\right) \right) , \\
&&\left( \left( e_{15},e_{24},e_{32}\right) ,\left(
U_{1},\{c_{3},c_{4}\},\{v_{1},v_{4}\}\right) \right) , \\
&&\left( \left( e_{14},e_{23},e_{33}\right) ,\left(
\{h_{3},h_{4},h_{5}\},\emptyset ,\{v_{2},v_{4}\}\right) \right) , \\
&&\left( \left( e_{12},e_{25},e_{32}\right) ,\left(
\{h_{1},h_{2}\},U_{2},\{v_{1},v_{4}\}\right) \right) , \\
&&\left( \left( e_{13},e_{21},e_{31}\right) ,\left( \emptyset
,\{c_{1},c_{2},c_{3}\},\{v_{2},v_{3}\}\right) \right) , \\
&&\left( \left( e_{11},e_{26},e_{32}\right) ,\left(
\{h_{3},h_{4},h_{5},h_{6}\},\{c_{4},c_{5}\},\{v_{1},v_{4}\}\right) \right) \}
\end{eqnarray*}
\end{example}

\begin{definition}
(\cite{Alkhazaleh}) For any soft multiset $(F_{A},E),$\ a pair $\left(
e_{ij},F_{ij}\right) $ is called a $U_{i}-$ soft multiset part for $\forall
e_{ij}\in a_{k},$ and $F_{ij}\subset F_{A}\left( A\right) $ is an
approximate value set, where $a_{k}\in A,$ $k=1,2,3,...,r,$ $%
i=1,2,...,m_{i}, $ $j=1,2,...,n_{j}.$
\end{definition}

\begin{example}
\label{ex3}Consider the soft multiset given in Example \ref{ex2}. Then,%
\begin{eqnarray*}
\left( e_{1j},F_{1j}\right) &=&\{\left(
e_{11},\{h_{3},h_{4},h_{5},h_{6}\}\right) ,\left(
e_{12},\{h_{1},h_{2}\}\right) ,\left( e_{13},\emptyset \right) , \\
&&\left( e_{14},\{h_{3},h_{4},h_{5}\}\right) ,\left( e_{15},U_{1}\right) \}
\end{eqnarray*}%
is a $U_{1}-$ soft multiset part of $(F_{A},E).$
\end{example}

\begin{definition}
(i) Let a pair $\left( e_{ij},F_{ij}\right) $\ be a $U_{i}-$ soft multiset
part of soft multiset $(F_{A},E),$ $A_{i}\subset E_{i}.$\ Then a subset of $%
U_{i}\times E_{i}$ is is called a relation form of $U_{i}-$ soft multiset
part of $(F_{A},E)$ which is uniquely defined by%
\begin{equation*}
R_{A_{i}}=\{\left( u_{ij},e_{ij}\right) :e_{ij}\in A_{i},u_{ij}\in
F_{ij}\left( e_{ij}\right) \}
\end{equation*}%
The characteristic function of $\chi _{R_{A_{i}}}$ is written by%
\begin{equation*}
\chi _{R_{A_{i}}}:U_{i}\times E_{i}\rightarrow \{0,1\},~~\chi
_{R_{A_{i}}}\left( u_{ij},e_{ij}\right) :=\left\{
\begin{array}{c}
1,\left( u_{ij},e_{ij}\right) \in R_{A_{i}} \\
0,\left( u_{ij},e_{ij}\right) \notin R_{A_{i}}%
\end{array}%
\right.
\end{equation*}%
(ii) If $U_{i}=\{u_{i1},u_{i2},...,u_{im_{i}}\},E_{i}=%
\{e_{i1},e_{i2},...,e_{in_{i}}\},$ then we call a matrix $\left[ a_{lk}^{i}%
\right] =\chi _{R_{A_{i}}}\left( u_{ik},e_{il}\right) ,1\leq k\leq
m_{i},1\leq l\leq n_{i}$ as an $m_{i}\times n_{i}$ soft matrix of $U_{i}-$
soft multiset part of $(F_{A},E).$
\end{definition}

\begin{example}
\label{ex4}Let us consider Example \ref{ex3}. Then%
\begin{eqnarray*}
R_{A_{1}} &=&\{\left( h_{3},e_{11}\right) ,\left( h_{4},e_{11}\right)
,\left( h_{5},e_{11}\right) ,\left( h_{6},e_{11}\right) ,\left(
h_{1},e_{12}\right) ,\left( h_{2},e_{12}\right) ,\left( h_{3},e_{14}\right) ,
\\
&&\left( h_{4},e_{14}\right) ,\left( h_{5},e_{14}\right) ,\left(
h_{1},e_{15}\right) ,\left( h_{2},e_{15}\right) ,\left( h_{3},e_{15}\right)
,\left( h_{4},e_{15}\right) ,\left( h_{5},e_{15}\right) ,\left(
h_{6},e_{15}\right) \}
\end{eqnarray*}%
Then $R_{A_{1}}$\ is presented by a table as in the following form:%
\begin{equation*}
\begin{tabular}{|l|l|l|l|l|l|}
\hline
$R_{A_{1}}$ & $e_{11}$ & $e_{12}$ & $e_{13}$ & $e_{14}$ & $e_{15}$ \\ \hline
$h_{1}$ & $\chi _{R_{A_{1}}}\left( h_{1},e_{11}\right) $ & $\chi
_{R_{A_{1}}}\left( h_{1},e_{12}\right) $ & $\chi _{R_{A_{1}}}\left(
h_{1},e_{13}\right) $ & $\chi _{R_{A_{1}}}\left( h_{1},e_{14}\right) $ & $%
\chi _{R_{A_{1}}}\left( h_{1},e_{15}\right) $ \\ \hline
$h_{2}$ & $\chi _{R_{A_{1}}}\left( h_{2},e_{11}\right) $ & $\chi
_{R_{A_{1}}}\left( h_{2},e_{12}\right) $ & $\chi _{R_{A_{1}}}\left(
h_{2},e_{13}\right) $ & $\chi _{R_{A_{1}}}\left( h_{2},e_{14}\right) $ & $%
\chi _{R_{A_{1}}}\left( h_{2},e_{15}\right) $ \\ \hline
$h_{3}$ & $\chi _{R_{A_{1}}}\left( h_{3},e_{11}\right) $ & $\chi
_{R_{A_{1}}}\left( h_{3},e_{12}\right) $ & $\chi _{R_{A_{1}}}\left(
h_{3},e_{13}\right) $ & $\chi _{R_{A_{1}}}\left( h_{3},e_{14}\right) $ & $%
\chi _{R_{A_{1}}}\left( h_{3},e_{15}\right) $ \\ \hline
$h_{4}$ & $\chi _{R_{A_{1}}}\left( h_{4},e_{11}\right) $ & $\chi
_{R_{A_{1}}}\left( h_{4},e_{12}\right) $ & $\chi _{R_{A_{1}}}\left(
h_{4},e_{13}\right) $ & $\chi _{R_{A_{1}}}\left( h_{4},e_{14}\right) $ & $%
\chi _{R_{A_{1}}}\left( h_{4},e_{15}\right) $ \\ \hline
$h_{5}$ & $\chi _{R_{A_{1}}}\left( h_{5},e_{11}\right) $ & $\chi
_{R_{A_{1}}}\left( h_{5},e_{12}\right) $ & $\chi _{R_{A_{1}}}\left(
h_{5},e_{13}\right) $ & $\chi _{R_{A_{1}}}\left( h_{5},e_{14}\right) $ & $%
\chi _{R_{A_{1}}}\left( h_{5},e_{15}\right) $ \\ \hline
$h_{6}$ & $\chi _{R_{A_{1}}}\left( h_{6},e_{11}\right) $ & $\chi
_{R_{A_{1}}}\left( h_{6},e_{12}\right) $ & $\chi _{R_{A_{1}}}\left(
h_{6},e_{13}\right) $ & $\chi _{R_{A_{1}}}\left( h_{6},e_{14}\right) $ & $%
\chi _{R_{A_{1}}}\left( h_{6},e_{15}\right) $ \\ \hline
\end{tabular}%
.
\end{equation*}%
Hence soft matrix $\left[ a_{lk}^{1}\right] $ of $U_{1}$ soft multiset part
of $(F_{A},E)$ is written by
\begin{equation*}
\left[ a_{lk}^{1}\right] _{6\times 5}=%
\begin{bmatrix}
0 & 1 & 0 & 0 & 1 \\
0 & 1 & 0 & 0 & 1 \\
1 & 0 & 0 & 1 & 1 \\
1 & 0 & 0 & 1 & 1 \\
1 & 0 & 0 & 1 & 1 \\
1 & 0 & 0 & 0 & 1%
\end{bmatrix}%
_{6\times 5}.
\end{equation*}
\end{example}

\begin{definition}
(\cite{Cagman3}) $\left[ a_{lk}\right] _{m\times n}$ is called a zero soft
matrix, denoted by $[0]$, if $a_{lk}=0$ for all $1\leq l\leq m$ and $1\leq
k\leq n.$
\end{definition}

\maketitle

\section{Soft matrices on soft multisets}

In this section, inspired by the above definitions to soft matrices and soft
multisets, first we will begin defining soft matrices on soft multisets and
its product and then we will give examples for these concepts.

\begin{definition}
\label{def1}Let $U_{i}=\{u_{i1},u_{i2},...,u_{im_{i}}\}$ be universes$%
,E_{i}=\{e_{i1},e_{i2},...,e_{in_{i}}\}$ be parameters for each $i\in I$, $%
R_{A_{i}}$ be a relation form of $U_{i}-$ soft multiset part of $(F_{A},E),$
$\left[ a_{lk}^{i}\right] _{m_{i}\times n_{i}},$ $1\leq k\leq m_{i},1\leq
l\leq n_{i}$ be soft matrix of $U_{i}-$ soft multiset part of $(F_{A},E).$
Then
\begin{equation*}
\left[ A_{lk}\right] _{m\times n}=\left[
\begin{tabular}{c|c|c|c}
$\left[ a_{lk}^{1}\right] _{m_{1}\times n_{1}}$ & $\left[ 0\right]
_{m_{1}\times n_{2}}$ & $\cdots $ & $\left[ 0\right] _{m_{1}\times n_{N}}$
\\ \hline
$\left[ 0\right] _{m_{2}\times n_{1}}$ & $\left[ a_{lk}^{2}\right]
_{m_{2}\times n_{2}}$ & $\cdots $ & $\left[ 0\right] _{m_{2}\times n_{N}}$
\\ \hline
$\cdots $ & $\cdots $ & $\ddots $ & $\cdots $ \\ \hline
$\left[ 0\right] _{m_{N}\times n_{1}}$ & $\left[ 0\right] _{m_{N}\times
n_{2}}$ & $\cdots $ & $\left[ a_{lk}^{N}\right] _{m_{N}\times n_{N}}$%
\end{tabular}%
\right] _{m\times n}
\end{equation*}%
is called a soft matrix of $(F_{A},E),$ where $m=m_{1}+m_{2}+\cdots
+m_{N},~n=n_{1}+n_{2}+\cdots +n_{N}.$
\end{definition}

\begin{example}
\label{ex6}Let us consider Example \ref{ex3}. \ The soft matrix $\left[
a_{lk}^{1}\right] $ of $U_{1}-$ soft multiset part of $(F_{A},E)$ is given
in Example \ref{ex4}. Similarly, we can obtain the soft matrices $\left[
a_{lk}^{2}\right] $ of $U_{2}-$ and $\left[ a_{lk}^{3}\right] $ of $U_{3}-$
soft multiset parts of $(F_{A},E)$ as shown below:%
\begin{equation*}
\left[ a_{lk}^{2}\right] _{5\times 6}=%
\begin{bmatrix}
1 & 0 & 0 & 0 & 1 & 0 \\
1 & 0 & 0 & 0 & 1 & 0 \\
1 & 0 & 0 & 1 & 1 & 0 \\
0 & 1 & 0 & 1 & 1 & 1 \\
0 & 1 & 0 & 0 & 1 & 1%
\end{bmatrix}%
_{5\times 6},\text{ }\left[ a_{lk}^{3}\right] _{4\times 5}=%
\begin{bmatrix}
0 & 1 & 0 & 1 & 0 \\
1 & 0 & 1 & 1 & 1 \\
1 & 0 & 0 & 0 & 1 \\
0 & 1 & 1 & 0 & 0%
\end{bmatrix}%
_{4\times 5}.
\end{equation*}%
Then soft matrix $\left[ A_{lk}\right] _{m\times n}$ of $(F_{A},E)$ can be
written%
\begin{equation*}
\left[ A_{lk}\right] _{15\times 16}=\left[
\begin{tabular}{c|c|c}
$%
\begin{bmatrix}
0 & 1 & 0 & 0 & 1 \\
0 & 1 & 0 & 0 & 1 \\
1 & 0 & 0 & 1 & 1 \\
1 & 0 & 0 & 1 & 1 \\
1 & 0 & 0 & 1 & 1 \\
1 & 0 & 0 & 0 & 1%
\end{bmatrix}%
_{6\times 5}$ & $%
\begin{bmatrix}
0 & 0 & 0 & 0 & 0 & 0 \\
0 & 0 & 0 & 0 & 0 & 0 \\
0 & 0 & 0 & 0 & 0 & 0 \\
0 & 0 & 0 & 0 & 0 & 0 \\
0 & 0 & 0 & 0 & 0 & 0 \\
0 & 0 & 0 & 0 & 0 & 0%
\end{bmatrix}%
_{6\times 6}$ & $%
\begin{bmatrix}
0 & 0 & 0 & 0 & 0 \\
0 & 0 & 0 & 0 & 0 \\
0 & 0 & 0 & 0 & 0 \\
0 & 0 & 0 & 0 & 0 \\
0 & 0 & 0 & 0 & 0 \\
0 & 0 & 0 & 0 & 0%
\end{bmatrix}%
_{6\times 5}$ \\ \hline
$%
\begin{bmatrix}
0 & 0 & 0 & 0 & 0 \\
0 & 0 & 0 & 0 & 0 \\
0 & 0 & 0 & 0 & 0 \\
0 & 0 & 0 & 0 & 0 \\
0 & 0 & 0 & 0 & 0%
\end{bmatrix}%
_{5\times 5}$ & $%
\begin{bmatrix}
1 & 0 & 0 & 0 & 1 & 0 \\
1 & 0 & 0 & 0 & 1 & 0 \\
1 & 0 & 0 & 1 & 1 & 0 \\
0 & 1 & 0 & 1 & 1 & 1 \\
0 & 1 & 0 & 0 & 1 & 1%
\end{bmatrix}%
_{5\times 6}$ & $%
\begin{bmatrix}
0 & 0 & 0 & 0 & 0 \\
0 & 0 & 0 & 0 & 0 \\
0 & 0 & 0 & 0 & 0 \\
0 & 0 & 0 & 0 & 0 \\
0 & 0 & 0 & 0 & 0%
\end{bmatrix}%
_{5\times 5}$ \\ \hline
$%
\begin{bmatrix}
0 & 0 & 0 & 0 & 0 \\
0 & 0 & 0 & 0 & 0 \\
0 & 0 & 0 & 0 & 0 \\
0 & 0 & 0 & 0 & 0%
\end{bmatrix}%
_{4\times 5}$ & $%
\begin{bmatrix}
0 & 0 & 0 & 0 & 0 & 0 \\
0 & 0 & 0 & 0 & 0 & 0 \\
0 & 0 & 0 & 0 & 0 & 0 \\
0 & 0 & 0 & 0 & 0 & 0%
\end{bmatrix}%
_{4\times 6}$ & $%
\begin{bmatrix}
0 & 1 & 0 & 1 & 0 \\
1 & 0 & 1 & 1 & 1 \\
1 & 0 & 0 & 0 & 1 \\
0 & 1 & 1 & 0 & 0%
\end{bmatrix}%
_{4\times 5}$%
\end{tabular}%
\right] _{15\times 16}
\end{equation*}
\end{example}
We make the following product definitions for soft matrices on soft
multisets, which are adapted from Definitions 7-10 in \cite{Cagman3}.

\begin{definition}
Let $\left[ a_{lk}^{i}\right] _{m_{i}\times n_{i}}$ be soft matrix of $%
U_{i}- $ soft multiset part of $(F_{A},E),$ $\left[ b_{lj}^{i}\right]
_{m_{i}\times n_{i}}$ be soft matrix of $U_{i}-$ soft multiset part of $%
(F_{B},E)$. Then \textbf{And }product of $\left[ a_{lk}^{i}\right]
_{m_{i}\times n_{i}}$ and $\left[ b_{lj}^{i}\right] _{m_{i}\times n_{i}}$ is
defined by%
\begin{equation*}
\left[ a_{lk}^{i}\right] _{m_{i}\times n_{i}}\wedge \left[ b_{lj}^{i}\right]
_{m_{i}\times n_{i}}=\left[ c_{lp}^{i}\right] _{m_{i}\times n_{i}^{2}}
\end{equation*}%
where $c_{lp}^{i}=\min \{a_{lk}^{i},b_{lj}^{i}\}$ such that $p=n_{i}\left(
k-1\right) +j.$
\end{definition}

\begin{definition}
Let $\left[ a_{lk}^{i}\right] _{m_{i}\times n_{i}}$ be soft matrix of $%
U_{i}- $ soft multiset part of $(F_{A},E),$ $\left[ b_{lj}^{i}\right]
_{m_{i}\times n_{i}}$ be soft matrix of $U_{i}-$ soft multiset part of $%
(F_{B},E)$. Then \textbf{Or }product of $\left[ a_{lk}^{i}\right]
_{m_{i}\times n_{i}}$ and $\left[ b_{lj}^{i}\right] _{m_{i}\times n_{i}}$ is
defined by%
\begin{equation*}
\left[ a_{lk}^{i}\right] _{m_{i}\times n_{i}}\vee \left[ b_{lj}^{i}\right]
_{m_{i}\times n_{i}}=\left[ c_{lp}^{i}\right] _{m_{i}\times n_{i}^{2}}
\end{equation*}%
where $c_{lp}^{i}=\max \{a_{lk}^{i},b_{lj}^{i}\}$ such that $p=n_{i}\left(
k-1\right) +j.$
\end{definition}

\begin{definition}
Let $\left[ a_{lk}^{i}\right] _{m_{i}\times n_{i}}$ be soft matrix of $%
U_{i}- $ soft multiset part of $(F_{A},E),$ $\left[ b_{lj}^{i}\right]
_{m_{i}\times n_{i}}$ be soft matrix of $U_{i}-$ soft multiset part of $%
(F_{B},E)$. Then \textbf{And-Not }product of $\left[ a_{lk}^{i}\right]
_{m_{i}\times n_{i}}$ and $\left[ b_{lj}^{i}\right] _{m_{i}\times n_{i}}$ is
defined by%
\begin{equation*}
\left[ a_{lk}^{i}\right] _{m_{i}\times n_{i}}\barwedge \left[ b_{lj}^{i}%
\right] _{m_{i}\times n_{i}}=\left[ c_{lp}^{i}\right] _{m_{i}\times
n_{i}^{2}}
\end{equation*}%
where $c_{lp}^{i}=\min \{a_{lk}^{i},1-b_{lj}^{i}\}$ such that $p=n_{i}\left(
k-1\right) +j.$
\end{definition}

\begin{definition}
Let $\left[ a_{lk}^{i}\right] _{m_{i}\times n_{i}}$ be soft matrix of $%
U_{i}- $ soft multiset part of $(F_{A},E),$ $\left[ b_{lj}^{i}\right]
_{m_{i}\times n_{i}}$ be soft matrix of $U_{i}-$ soft multiset part of $%
(F_{B},E)$. Then \textbf{Or-Not }product of $\left[ a_{lk}^{i}\right]
_{m_{i}\times n_{i}}$ and $\left[ b_{lj}^{i}\right] _{m_{i}\times n_{i}}$ is
defined by%
\begin{equation*}
\left[ a_{lk}^{i}\right] _{m_{i}\times n_{i}}\veebar \left[ b_{lj}^{i}\right]
_{m_{i}\times n_{i}}=\left[ c_{lp}^{i}\right] _{m_{i}\times n_{i}^{2}}
\end{equation*}%
where $c_{lp}^{i}=\max \{a_{lk}^{i},1-b_{lj}^{i}\}$ such that $p=n_{i}\left(
k-1\right) +j.$
\end{definition}

\begin{example}
\label{ex5}As in Example \ref{ex2}, Mr. X is joining with Mrs. X for
accommodation purchase, transportation purchase, and venue to hold a wedding
celebration respectively. The set of choice parameters for him is%
\begin{eqnarray*}
B &=&\{b_{1}=\left( e_{11},e_{25},e_{31}\right) ,b_{2}=\left(
e_{13},e_{22},e_{33}\right) ,b_{3}=\left( e_{12},e_{26},e_{32}\right)
,b_{4}=\left( e_{14},e_{24},e_{34}\right) , \\
b_{5} &=&\left( e_{15},e_{23},e_{35}\right) ,b_{6}=\left(
e_{11},e_{21},e_{31}\right) \}\subset E.
\end{eqnarray*}%
Then the soft multiset $(F_{B},E)$ is given by%
\begin{eqnarray*}
(F_{B},E) &=&\{\left( \left( e_{11},e_{25},e_{31}\right) ,\left(
U_{1},\{c_{4},c_{5}\},\{v_{1},v_{2},v_{3}\}\right) \right) , \\
&&\left( \left( e_{13},e_{22},e_{33}\right) ,\left(
\{h_{2},h_{3},h_{4},h_{5}\},\{c_{1},c_{2}\},\{v_{2},v_{4}\}\right) \right) ,
\\
&&\left( \left( e_{12},e_{26},e_{32}\right) ,\left( \emptyset
,\{c_{4},c_{5}\},\{v_{4}\}\right) \right) , \\
&&\left( \left( e_{14},e_{24},e_{34}\right) ,\left(
\{h_{1},h_{2},h_{3}\},\{c_{4},c_{5}\},\{v_{2},v_{3}\}\right) \right) , \\
&&\left( \left( e_{15},e_{23},e_{35}\right) ,\left(
\{h_{1},h_{2},h_{5},h_{6}\},\{c_{1},c_{2},c_{3},c_{4}\},\{v_{2}\}\right)
\right) , \\
&&\left( \left( e_{11},e_{21},e_{31}\right) ,\left(
U_{1},\{c_{3},c_{4},c_{5}\},\{v_{1},v_{2},v_{3}\}\right) \right) \},
\end{eqnarray*}%
and the soft matrix $\left[ b_{lk}^{1}\right] $ of $U_{1}-$ soft multiset
part of $(F_{B},E)$ is obtained as shown below:%
\begin{equation*}
\left[ b_{lk}^{1}\right] _{6\times 5}=%
\begin{bmatrix}
1 & 0 & 0 & 1 & 1 \\
1 & 0 & 1 & 1 & 1 \\
1 & 0 & 1 & 1 & 0 \\
1 & 0 & 1 & 0 & 0 \\
1 & 0 & 1 & 0 & 1 \\
1 & 0 & 0 & 0 & 1%
\end{bmatrix}%
_{6\times 5}.
\end{equation*}%
Then we have \textbf{And }product of $\left[ a_{ij}^{1}\right] _{m\times n}$
of $U_{1}-$ soft multiset part of $(F_{A},E)$ and $\left[ b_{ik}^{1}\right]
_{m\times n}$ of $U_{1}-$ soft multiset part of $(F_{B},E)$ in the following:%

%\begin{figure}[H]
%  % Requires \usepackage{graphicx}
%  \includegraphics[width=\textwidth]{fig_SoftSetC1And.ps}\\
% % \caption{}\label{}
%\end{figure}

\begin{figure}
\centering
\includegraphics[width=\textwidth]{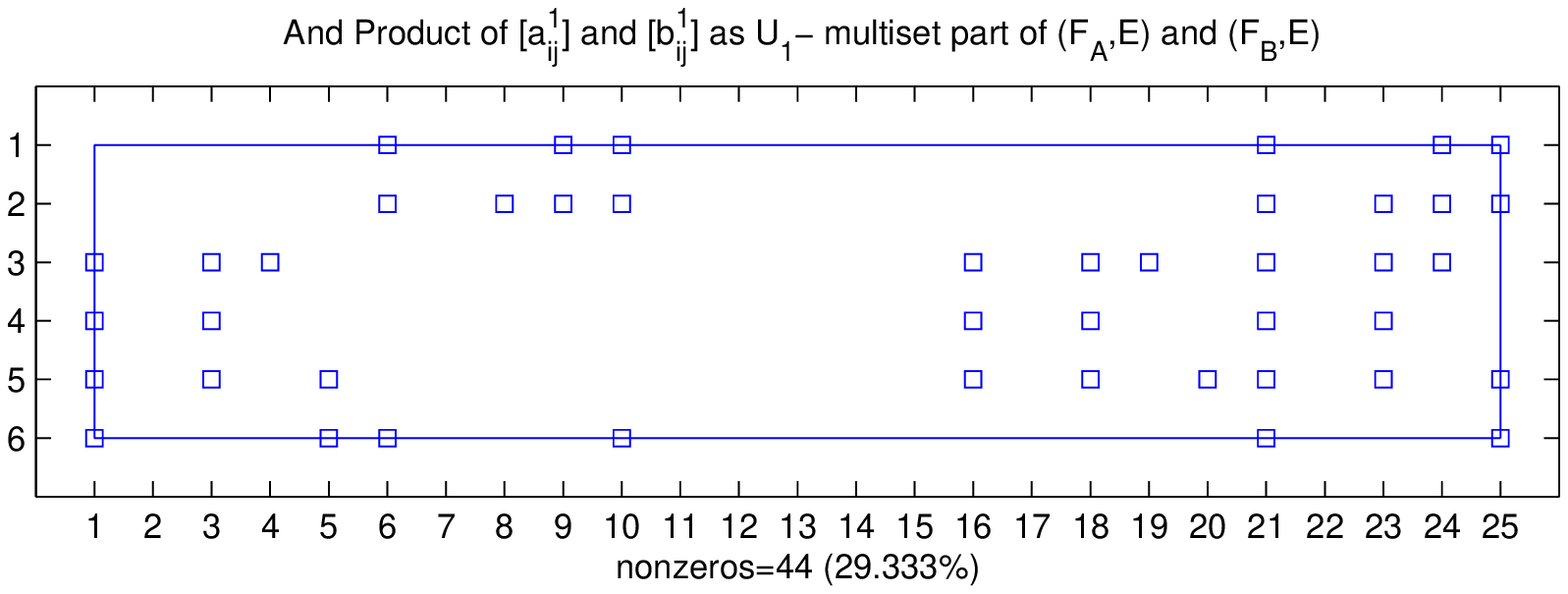}
\end{figure}
\begin{figure}[H]
\centering
\includegraphics[width=\textwidth]{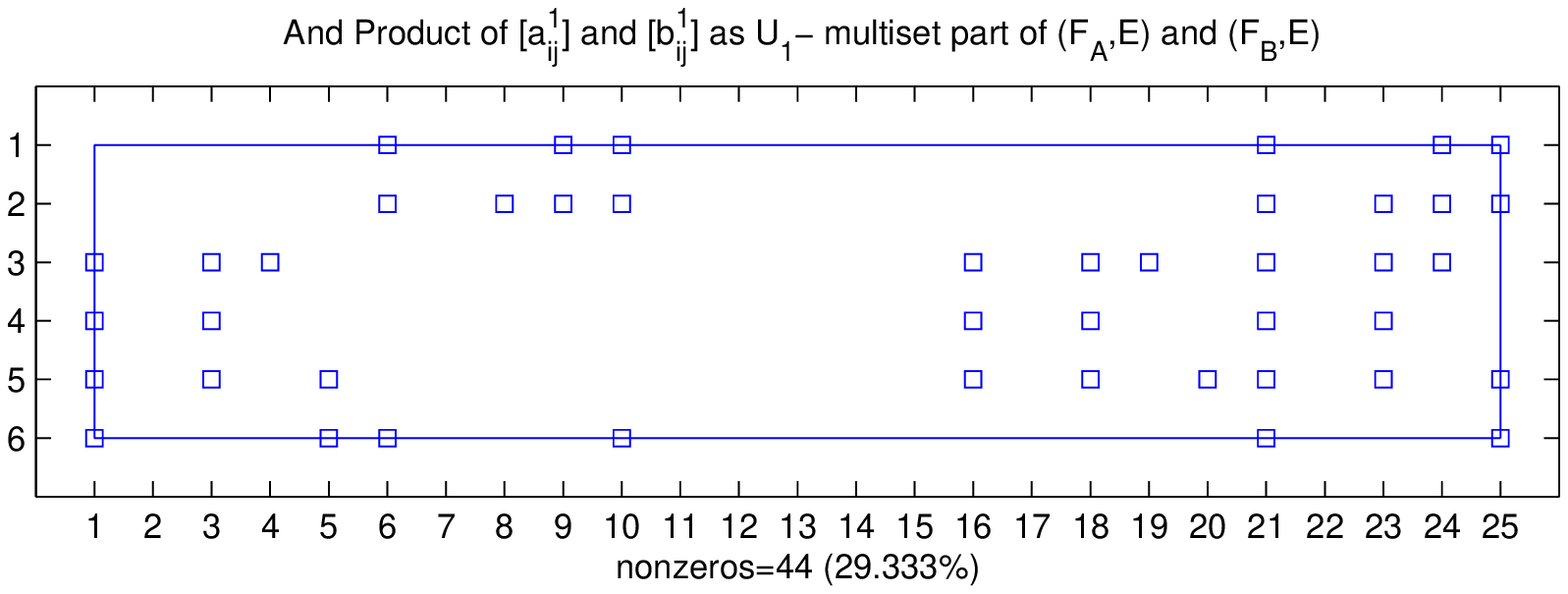}
\end{figure}
\textbf{Or }product of $%
\left[ a_{ij}^{1}\right] _{m\times n}$ of $U_{1}$ soft multiset part of $%
(F_{A},E)$ and $\left[ b_{ik}^{1}\right] _{m\times n}$ of $U_{1}$ soft
multiset part of $(F_{B},E)$ in the following:
\begin{figure}[H]
\centering\includegraphics[width=\textwidth]{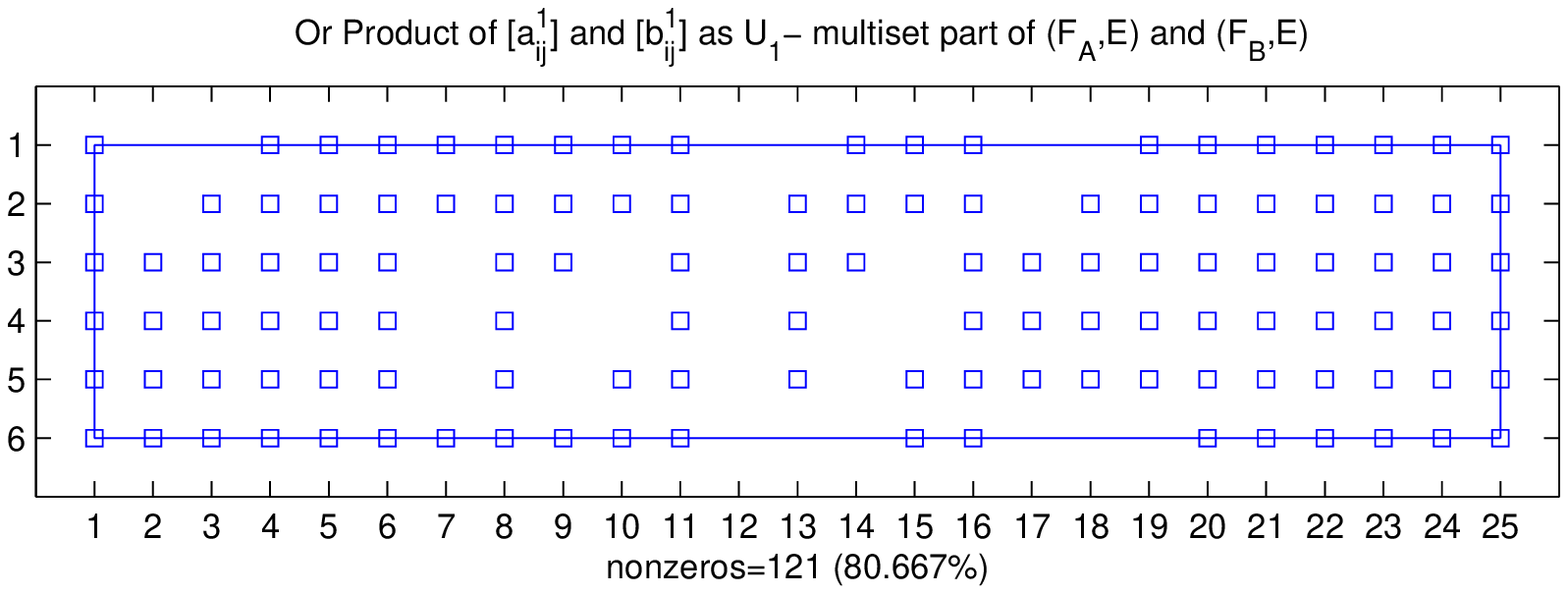}
\end{figure}
\textbf{And-Not }product of $\left[ a_{ij}^{1}\right] _{m\times n}$ of $U_{1}
$ soft multiset part of $(F_{A},E)$ and $\left[ b_{ik}^{1}\right] _{m\times
n}$ of $U_{1}$ soft multiset part of $(F_{B},E)$ in the following:
\begin{figure}[H]
\centering
\includegraphics[width=\textwidth]{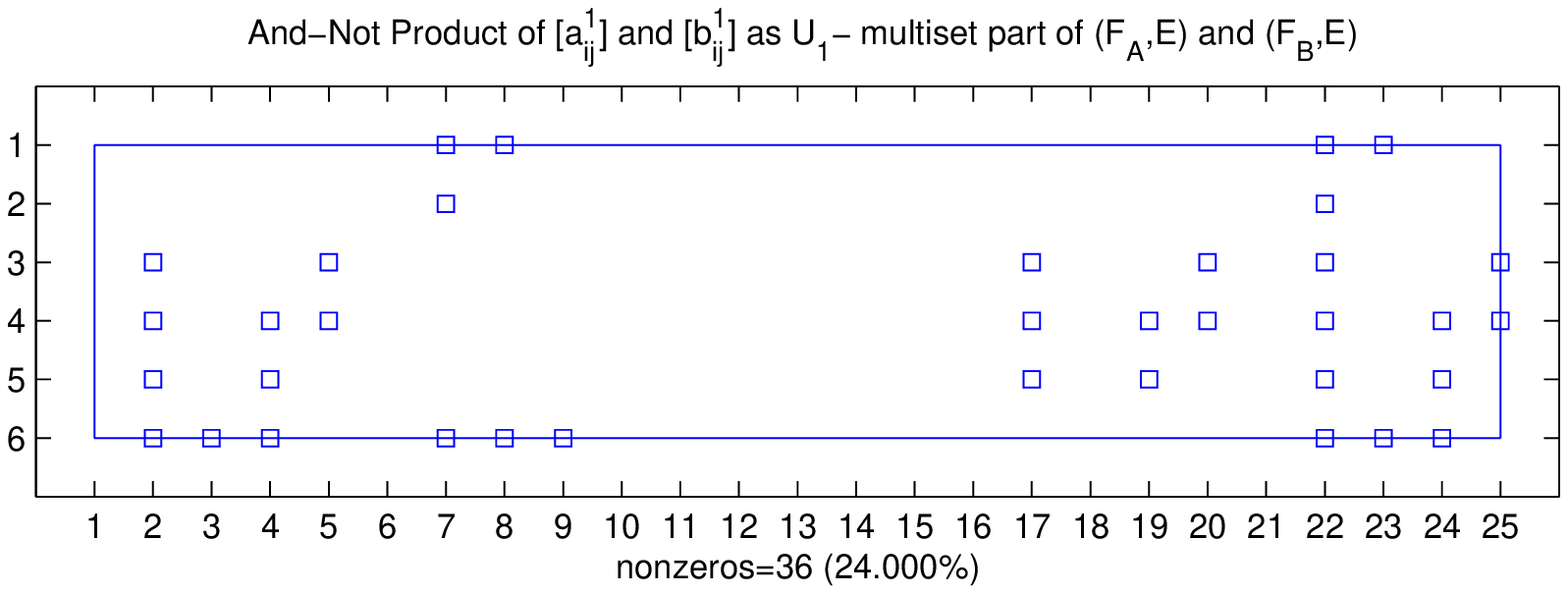}
\end{figure}
\textbf{Or-Not }product
of $\left[ a_{ij}^{1}\right] _{m\times n}$ of $U_{1}$ soft multiset part of $%
(F_{A},E)$ and $\left[ b_{ik}^{1}\right] _{m\times n}$ of $U_{1}$ soft
multiset part of $(F_{B},E)$ in the following:
\begin{figure}[H]
\centering
\includegraphics[width=\textwidth]{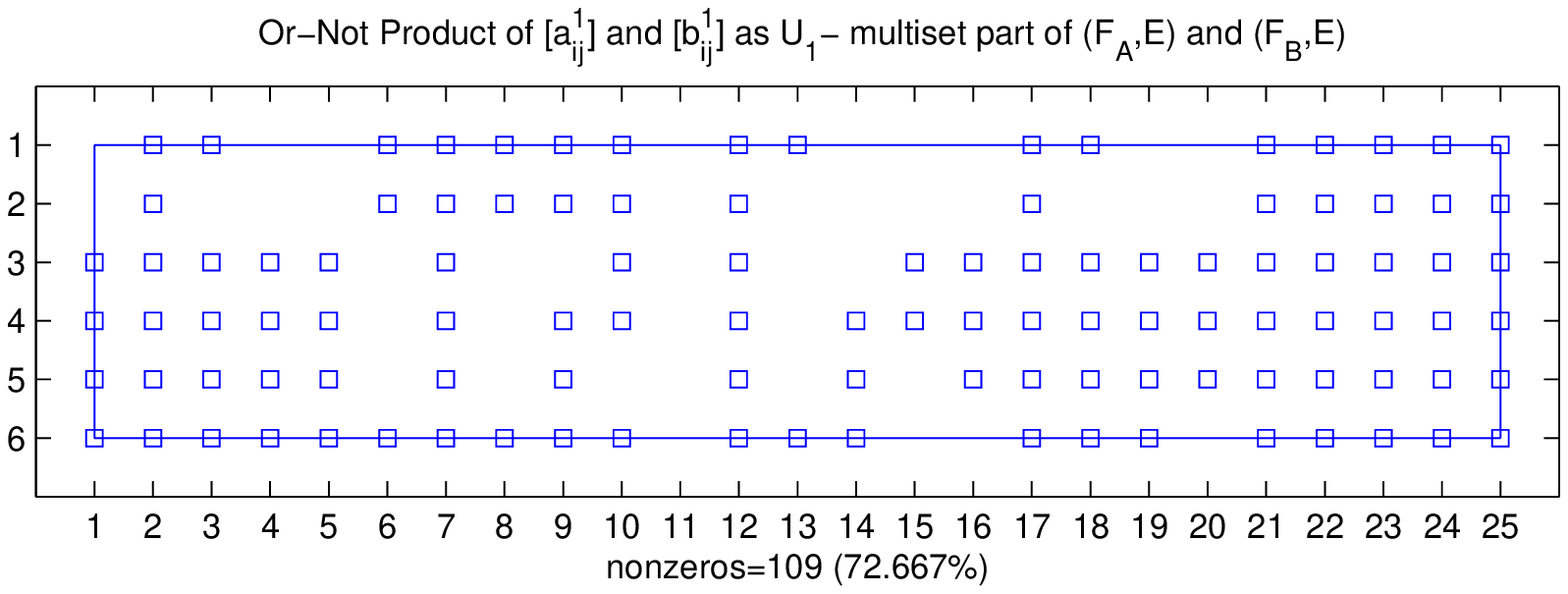}
\end{figure}
Here, squares show the value $1$ for the product of $\left[ a_{ij}^{1}\right] _{m\times n}$
of $U_{1}-$ soft multiset part of $(F_{A},E)$ and $\left[ b_{ik}^{1}\right]
_{m\times n}$ of $U_{1}-$ soft multiset part of $(F_{B},E)$.

\end{example}

\begin{definition}
\label{def2}Let $\left[ A_{lk}\right] _{m\times n},\left[ B_{lk}\right]
_{m\times n}$ be the soft matrices of soft multisets$\ (F_{A},E)$ and $%
\left( F_{B},E\right) ,$ respectively. Then \textbf{And }product of $\left[
A_{lk}\right] _{m\times n}$ and $\left[ B_{lk}\right] _{m\times n}$ is
defined by%
\begin{equation*}
\left[ A_{lk}\right] _{m\times n}\wedge \left[ B_{lk}\right] _{m\times n}=%
\left[
\begin{tabular}{c|c|c|c}
$\left[ a_{lk}^{1}\right] _{m_{1}\times n_{1}}\wedge \left[ b_{lk}^{1}\right]
_{m_{1}\times n_{1}}$ & $\left[ 0\right] _{m_{1}\times n_{2}^{2}}$ & $\cdots
$ & $\left[ 0\right] _{m_{1}\times n_{N}^{2}}$ \\ \hline
$\left[ 0\right] _{m_{2}\times n_{1}^{2}}$ & $\left[ a_{lk}^{2}\right]
_{m_{2}\times n_{2}}\wedge \left[ b_{lk}^{2}\right] _{m_{2}\times n_{2}}$ & $%
\cdots $ & $\left[ 0\right] _{m_{2}\times n_{N}^{2}}$ \\ \hline
$\cdots $ & $\cdots $ & $\ddots $ & $\cdots $ \\ \hline
$\left[ 0\right] _{m_{N}\times n_{1}^{2}}$ & $\left[ 0\right] _{m_{N}\times
n_{2}^{2}}$ & $\cdots $ & $\left[ a_{lk}^{N}\right] _{m_{N}\times
n_{N}}\wedge \left[ b_{lk}^{N}\right] _{m_{N}\times n_{N}}$%
\end{tabular}%
\right] _{m\times ns},
\end{equation*}%
\textbf{Or }product of $\left[ A_{lk}\right] _{m\times n}$ and $\left[ B_{lk}%
\right] _{m\times n}$ is defined by%
\begin{equation*}
\left[ A_{lk}\right] _{m\times n}\vee \left[ B_{lk}\right] _{m\times n}=%
\left[
\begin{tabular}{c|c|c|c}
$\left[ a_{lk}^{1}\right] _{m_{1}\times n_{1}}\vee \left[ b_{lk}^{1}\right]
_{m_{1}\times n_{1}}$ & $\left[ 0\right] _{m_{1}\times n_{2}^{2}}$ & $\cdots
$ & $\left[ 0\right] _{m_{1}\times n_{N}^{2}}$ \\ \hline
$\left[ 0\right] _{m_{2}\times n_{1}^{2}}$ & $\left[ a_{lk}^{2}\right]
_{m_{2}\times n_{2}}\vee \left[ b_{lk}^{2}\right] _{m_{2}\times n_{2}}$ & $%
\cdots $ & $\left[ 0\right] _{m_{2}\times n_{N}^{2}}$ \\ \hline
$\cdots $ & $\cdots $ & $\ddots $ & $\cdots $ \\ \hline
$\left[ 0\right] _{m_{N}\times n_{1}^{2}}$ & $\left[ 0\right] _{m_{N}\times
n_{2}^{2}}$ & $\cdots $ & $\left[ a_{lk}^{N}\right] _{m_{N}\times n_{N}}\vee %
\left[ b_{lk}^{N}\right] _{m_{N}\times n_{N}}$%
\end{tabular}%
\right] _{m\times ns},
\end{equation*}%
\textbf{And-Not }product of $\left[ A_{lk}\right] _{m\times n}$ and $\left[
B_{lk}\right] _{m\times n}$ is defined by%
\begin{equation*}
\left[ A_{lk}\right] _{m\times n}\barwedge \left[ B_{lk}\right] _{m\times n}=%
\left[
\begin{tabular}{c|c|c|c}
$\left[ a_{lk}^{1}\right] _{m_{1}\times n_{1}}\barwedge \left[ b_{lk}^{1}%
\right] _{m_{1}\times n_{1}}$ & $\left[ 0\right] _{m_{1}\times n_{2}^{2}}$ &
$\cdots $ & $\left[ 0\right] _{m_{1}\times n_{N}^{2}}$ \\ \hline
$\left[ 0\right] _{m_{2}\times n_{1}^{2}}$ & $\left[ a_{lk}^{2}\right]
_{m_{2}\times n_{2}}\barwedge \left[ b_{lk}^{2}\right] _{m_{2}\times n_{2}}$
& $\cdots $ & $\left[ 0\right] _{m_{2}\times n_{N}^{2}}$ \\ \hline
$\cdots $ & $\cdots $ & $\ddots $ & $\cdots $ \\ \hline
$\left[ 0\right] _{m_{N}\times n_{1}^{2}}$ & $\left[ 0\right] _{m_{N}\times
n_{2}^{2}}$ & $\cdots $ & $\left[ a_{lk}^{N}\right] _{m_{N}\times
n_{N}}\barwedge \left[ b_{lk}^{N}\right] _{m_{N}\times n_{N}}$%
\end{tabular}%
\right] _{m\times ns},
\end{equation*}%
\textbf{Or-Not }product of $\left[ A_{lk}\right] _{m\times n}$ and $\left[
B_{lk}\right] _{m\times n}$ is defined by%
\begin{equation*}
\left[ A_{lk}\right] _{m\times n}\veebar \left[ B_{lk}\right] _{m\times n}=%
\left[
\begin{tabular}{c|c|c|c}
$\left[ a_{lk}^{1}\right] _{m_{1}\times n_{1}}\veebar \left[ b_{lk}^{1}%
\right] _{m_{1}\times n_{1}}$ & $\left[ 0\right] _{m_{1}\times n_{2}^{2}}$ &
$\cdots $ & $\left[ 0\right] _{m_{1}\times n_{N}^{2}}$ \\ \hline
$\left[ 0\right] _{m_{2}\times n_{1}^{2}}$ & $\left[ a_{lk}^{2}\right]
_{m_{2}\times n_{2}}\veebar \left[ b_{lk}^{2}\right] _{m_{2}\times n_{2}}$ &
$\cdots $ & $\left[ 0\right] _{m_{2}\times n_{N}^{2}}$ \\ \hline
$\cdots $ & $\cdots $ & $\ddots $ & $\cdots $ \\ \hline
$\left[ 0\right] _{m_{N}\times n_{1}^{2}}$ & $\left[ 0\right] _{m_{N}\times
n_{2}^{2}}$ & $\cdots $ & $\left[ a_{lk}^{N}\right] _{m_{N}\times
n_{N}}\veebar \left[ b_{lk}^{N}\right] _{m_{N}\times n_{N}}$%
\end{tabular}%
\right] _{m\times ns},
\end{equation*}%
where $m=m_{1}+m_{2}+\cdots +m_{N},$ $ns=n_{1}^{2}+n_{2}^{2}+\cdots
+n_{N}^{2}$.
\end{definition}

\begin{example}
\label{ex7}\label{ex7}Let us consider Example \ref{ex6} and \ref{ex5}. \textbf{And }%
product of $\left[ A_{lk}\right] _{15\times 16}$ and $\left[ B_{lk}\right]
_{15\times 16}$ is $15\times 86$ soft matrix whose visualize sparsity
pattern is given in the following:
\begin{figure}[H]
\centering
\includegraphics[width=\textwidth]{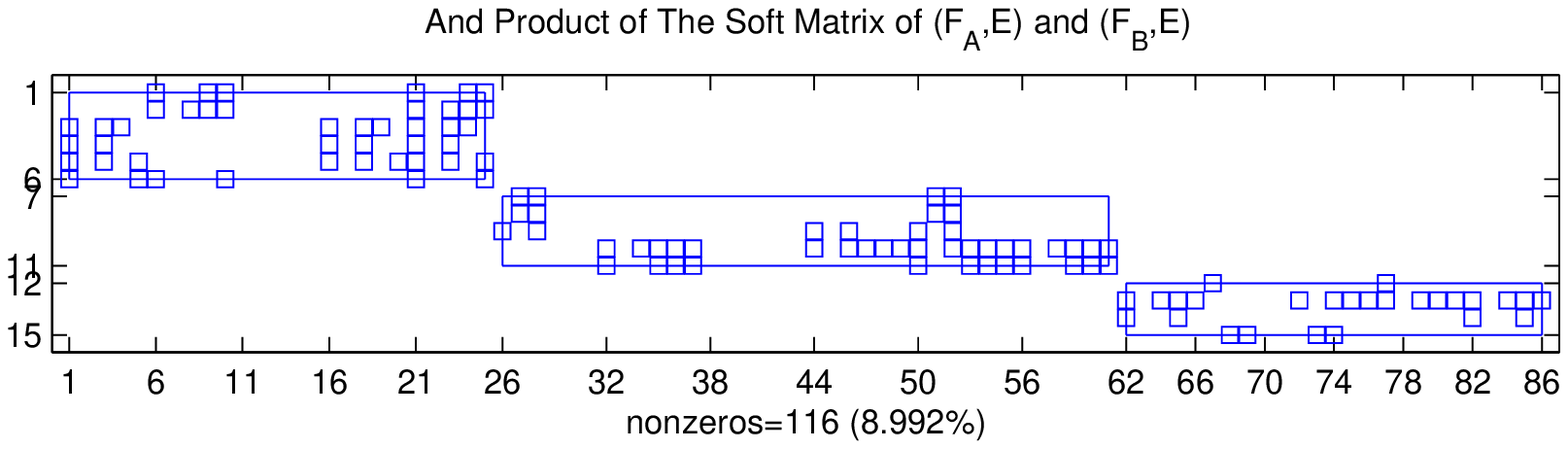}
\end{figure}
\textbf{Or }product of $\left[
A_{lk}\right] _{15\times 16}$ and $\left[ B_{lk}\right] _{15\times 16}$ is $%
15\times 86$ soft matrix whose visualize sparsity pattern is given in the
following:
\begin{figure}[H]
\centering
\includegraphics[width=\textwidth]{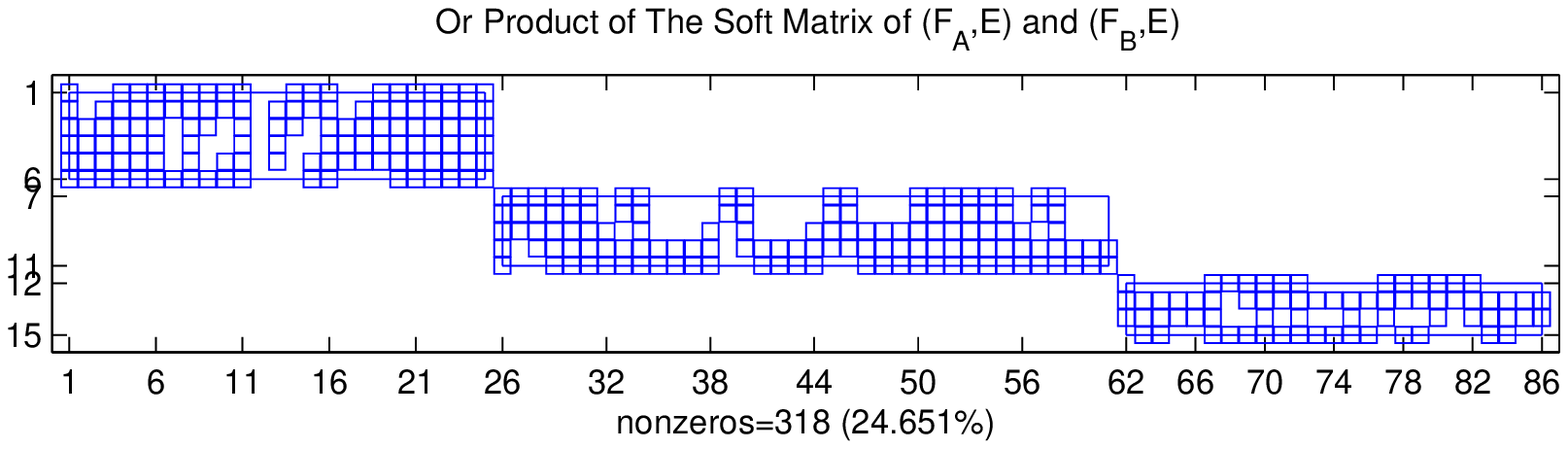}
\end{figure}
\textbf{And-Not }product of $%
\left[ A_{lk}\right] _{15\times 16}$ and $\left[ B_{lk}\right] _{15\times
16} $ is $15\times 86$ soft matrix whose visualize sparsity pattern is given
in the following:
\begin{figure}[H]
\centering
\includegraphics[width=\textwidth]{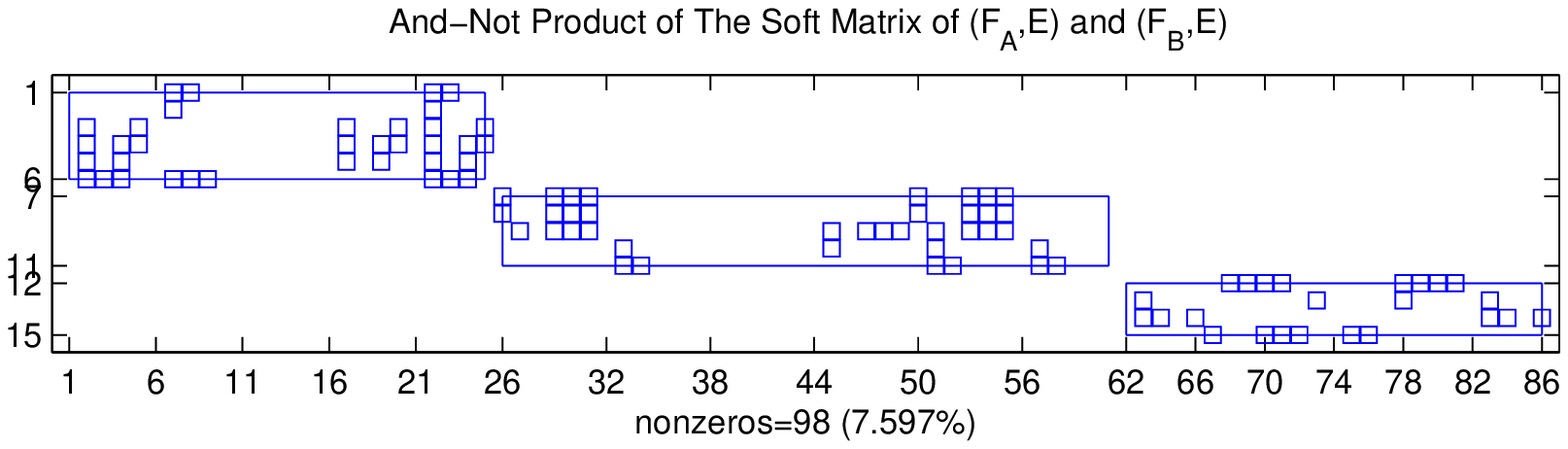}
\end{figure}
\textbf{Or-Not }product of
$\left[ A_{lk}\right] _{15\times 16}$ and $\left[ B_{lk}\right] _{15\times
16}$ is $15\times 86$ soft matrix whose visualize sparsity pattern is given
in the following:
\begin{figure}[H]
\centering
\includegraphics[width=\textwidth]{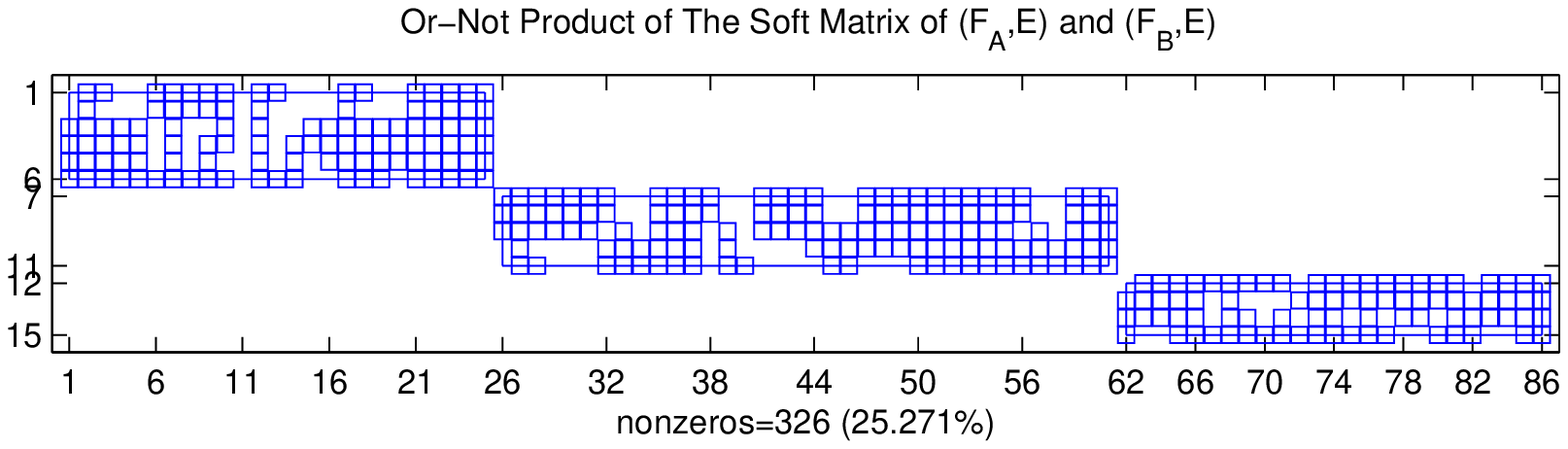}
\end{figure}
Here, squares show the
value $1$ for the product of $\left[ A_{lk}\right] _{m\times n}$ and $\left[ B_{lk}\right]
_{m\times n}$ as the soft matrices of soft multisets$\ (F_{A},E)$ and $%
\left( F_{B},E\right) ,$ respectively.
\end{example}

\section{Application}

Now we use the algorithm to solve our original problem.

\begin{example}
Suppose that a married couple, Mr. X and Mrs. X, are considering for
accommodation purchase, transportation purchase, and venue to hold a wedding
celebration respectively. The set of choice parameters for them are given in
Example \ref{ex2} and Example \ref{ex5}. We now select the house(s), car(s)
and hotel(s) on the sets of partners' parameters by using the above
algorithm as follows:

\begin{description}
\item[Step 1] First, Mrs. X and Mr. X have to choose the sets of their
parameters, given in Example \ref{ex2} and Example \ref{ex5}, respectively;

\item[Step 2] Then we can write the soft multisets $\left( F_{A},E\right)
,\left( F_{B},E\right) $, given in Example \ref{ex2} and Example \ref{ex5},
respectively;

\item[Step 3] Next we find the soft matrix $\left[ a_{lk}^{i}\right] $ of $%
U_{i}$ soft multiset part of $(F_{A},E)$ and the soft matrix $\left[
b_{lk}^{i}\right] $ of $U_{i}$ soft multiset part of $(F_{B},E)$, given in
Example \ref{ex6} and Example \ref{ex5}, respectively;

\item[Step 4] By using Definition \ref{def1}, we construct the soft matrix $%
\left[ A_{lk}\right] _{m\times n}$ of soft set $(F_{A},E)$ and the soft
matrix $\left[ B_{lk}\right] _{m\times n}$ of soft set $(F_{A},E);$
\begin{figure}[H]
\centering
\includegraphics[width=0.5\textwidth]{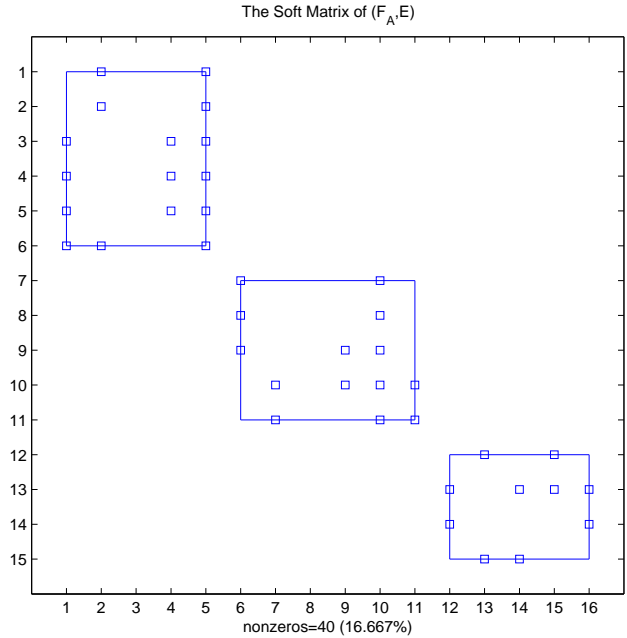}
\caption{Visualize Sparsity Pattern of the Soft
Matrix of $\left( F_{A},E\right) $}
\end{figure}

\begin{figure}[H]
\centering
\includegraphics[width=0.5\textwidth]{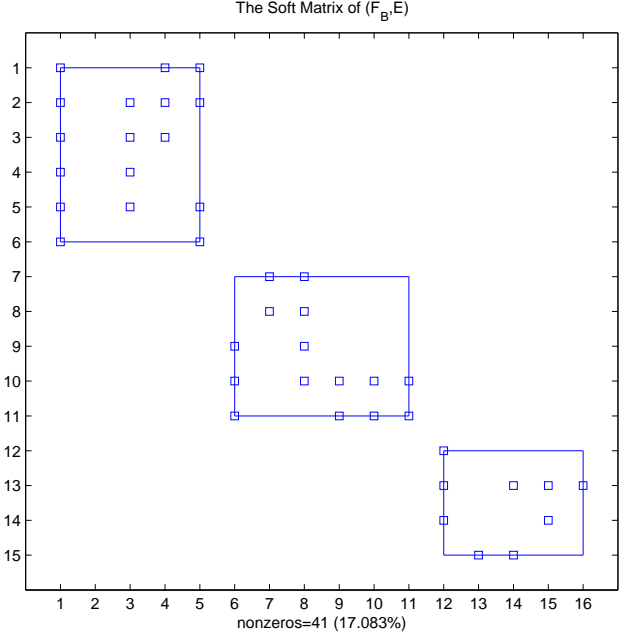}
\caption{Visualize Sparsity Pattern of the Soft
Matrix of $\left( F_{B},E\right) $}
\end{figure}

\item[Step 5] We construct \textbf{And }product of $\left[ A_{lk}\right]
_{m\times n}$ and $\left[ B_{lk}\right] _{m\times n}$ by using Definition %
\ref{def2}, denote by $\left[ C_{lj}\right] _{m\times ns},$ given by Example %
\ref{ex7};

\item[Step 6] We find the sets $I_{k}^{\left( i\right) }$ as%
\begin{eqnarray*}
I_{1}^{\left( 1\right) } &=&\{1;3;4;5\},I_{2}^{\left( 1\right)
}=\{6;8;9;10\},I_{3}^{\left( 1\right) }=\emptyset ,I_{4}^{\left( 1\right)
}=\{16;18;19;20\},I_{5}^{\left( 1\right) }=\{21;23;24;25\} \\
I_{1}^{\left( 2\right) } &=&\{26;27;28\},I_{2}^{\left( 2\right)
}=\{32;34;35;36;37\},I_{3}^{\left( 2\right) }=\emptyset ,I_{4}^{\left(
2\right) }=\{44;46;47;48;49\}, \\
I_{5}^{\left( 2\right) } &=&\{50;51;52;53;54;55\},I_{6}^{\left( 2\right)
}=\{56;58;59;60;61\} \\
I_{1}^{\left( 3\right) } &=&\{62;64;65;66\},I_{2}^{\left( 3\right)
}=\{67;68;69\},I_{3}^{\left( 3\right) }=\{72;73;74;75;76\}, \\
I_{4}^{\left( 3\right) } &=&\{77;79;80;81\},I_{5}^{\left( 3\right)
}=\{82;84;85;86\}
\end{eqnarray*}%
and construct decision function%
\begin{eqnarray*}
\left[ w_{l,k}^{\left( 1\right) }\right]  &=&%
\begin{bmatrix}
0 & 0 & 0 & 0 & 0 \\
0 & 1 & 0 & 0 & 1 \\
0 & 0 & 0 & 0 & 0 \\
0 & 0 & 0 & 0 & 0 \\
0 & 0 & 0 & 0 & 0 \\
0 & 0 & 0 & 0 & 0%
\end{bmatrix}%
_{6\times 5},\left[ w_{l,k}^{\left( 2\right) }\right] =%
\begin{bmatrix}
0 & 0 & 0 & 0 & 0 & 0 \\
0 & 0 & 0 & 0 & 0 & 0 \\
0 & 0 & 0 & 0 & 0 & 0 \\
0 & 1 & 0 & 1 & 0 & 1 \\
0 & 0 & 0 & 0 & 0 & 0%
\end{bmatrix}%
_{5\times 6},\left[ w_{l,k}^{\left( 3\right) }\right] =%
\begin{bmatrix}
0 & 0 & 0 & 0 & 0 \\
1 & 0 & 0 & 1 & 1 \\
0 & 0 & 0 & 0 & 0 \\
0 & 0 & 0 & 0 & 0%
\end{bmatrix}%
_{4\times 5} \\
\left[ v_{l}^{\left( 1\right) }\right]  &=&%
\begin{bmatrix}
0 \\
\mathbf{1} \\
0 \\
0 \\
0 \\
0%
\end{bmatrix}%
\begin{array}{c}
\rightarrow h_{1} \\
\rightarrow h_{2} \\
\rightarrow h_{3} \\
\rightarrow h_{4} \\
\rightarrow h_{5} \\
\rightarrow h_{6}%
\end{array}%
,\left[ v_{l}^{\left( 2\right) }\right] =%
\begin{bmatrix}
0 \\
0 \\
0 \\
\mathbf{1} \\
0%
\end{bmatrix}%
\begin{array}{c}
\rightarrow c_{1} \\
\rightarrow c_{2} \\
\rightarrow c_{3} \\
\rightarrow c_{4} \\
\rightarrow c_{5}%
\end{array}%
,\left[ v_{l}^{\left( 3\right) }\right] =%
\begin{bmatrix}
0 \\
\mathbf{1} \\
0 \\
0%
\end{bmatrix}%
\begin{array}{c}
\rightarrow v_{1} \\
\rightarrow v_{2} \\
\rightarrow v_{3} \\
\rightarrow v_{4}%
\end{array}%
\end{eqnarray*}

\item[Step 7] We find the optimum set of $U_{1}=\{h_{2}\},U_{2}=\{c_{4}%
\},U_{3}=\{v_{2}\}$
\end{description}
\end{example}


\begin{thebibliography}{20}
\bibitem{Ali}  M.I. Ali, F. Feng, X. Liu, W.K. Min, M. Shabir,  \emph{On Some new operations in Soft Set theory},
Computers and Mathematics with Applications, 57(2009), 1547-1553.

\bibitem{Alkhazaleh}  S. Alkhazaleh, A.R. Salleh,   \emph{Soft Multisets Theory},
Applied Mathematical Sciences, 5/72(2011), 3561 - 3573.


\bibitem{Alkhazaleh2}  S. Alkhazaleh, A.R. Salleh,   \emph{Fuzzy Soft Multiset Theory},
Abstract and Applied Analysis, 2012/350603(2012), 20 pages.

\bibitem{Basu}  T.M. Basu, N.K. Mahapatra, S.K. Mondal,   \emph{Different Types of Matrices in Fuzzy Soft Set Theory and Their Application in Decision
Making Problems},
IRACST -- Engineering Science and Technology: An International Journal (ESTIJ), 2/3(2012), 389-398.


\bibitem{Cagman} N. \c{C}a\u{g}man,  S.  Engino\u{g}lu ,   \emph{Soft set theory and uni-int decision making},
European Journal of Operational Research, 207(2010), 848-855.

\bibitem{Cagman2} N. \c{C}a\u{g}man,  S.  Engino\u{g}lu ,   \emph{Fuzzy Soft Matrix Theory and Its Application in Decision Making},
Iranian Journal of Fuzzy Systems, 9/1(2012), 109-119.


\bibitem{Cagman3} N. \c{C}a\u{g}man,  S.  Engino\u{g}lu ,   \emph{Soft matrix theory and its decision making},
Computers and Mathematics with Applications, 59(2010), 3308-3314.

\bibitem{Feng} F. Feng, X.  Liu,  V. L.  Fotea, Y. B. Jun ,   \emph{Soft sets and soft rough sets},
Information Sciences, 181(2011), 1125-1137.

\bibitem{Caras} C.A. Gunduz, S. Bayramov,   \emph{Some Results on Fuzzy Soft Topological Spaces},
Mathematical Problems in Engineering, 2013(2013), 10 pages.

\bibitem{Caras1} C.A. Gunduz, A. Sonmez,  H. Cakalli,   \emph{On Soft Mappings},
arXiv:1305.4545(2013), 12 pages.

\bibitem{Han} B. H. Han, S.L.  Geng,   \emph{Pruning method for optimal solutions of int m-int n decision making scheme},
European Journal of Operational Research, 231(2013), 779-783.

\bibitem{Herawan} T. Herawan, R.  Ghazali, M.M. Deris,    \emph{Soft Set Theoretic Approach for Dimensionality Reduction},
International Journal of Database Theory and Application, 3/2(2010), 47-60.

\bibitem{Li} Z. Li, T. Xie ,    \emph{The relationship among soft sets, soft rough sets and topologies},
Soft Comput., 10.1007/s00500-013-1108-5(2013).

\bibitem{Maji} P.K. Maji, R. Biswas, A.R. Roy,     \emph{Soft Set Theory},
Comput. Math. Appl., 45(2003), 555-562.

\bibitem{Maji2} P.K. Maji, A.R. Roy,     \emph{An Application of Soft Sets in A Decision Making Problem},
Computers and Mathematics with Applications, 44(2002), 1077-1083.

\bibitem{Mamat} R. Mamat, T. Herawan,  M.M. Deris,     \emph{MAR: Maximum Attribute Relative of soft set for clustering attribute selection},
Knowledge-Based Systems, 52(2013), 11-20.

\bibitem{Molodtsov} D. Molodtsov,     \emph{Soft set theory first result},
Computers and Mathematics with Applications, 37(1999), 19-31.

\bibitem{Mondal2} J.I. Mondal, T.K.  Roy,      \emph{Theory of Fuzzy Soft Matrix and its Multi Criteria in Decision Making Based on Three Basic t-Norm Operators},
International Journal of Innovative Research in Science, Engineering and Technology, 2/10(2013), 5715- 5723.

\bibitem{Mondal} S. Mondal,  M. Pal,     \emph{Soft matrices},
African Journal of Mathematics and Computer Science Research, 4/13(2011), 379-388.

\bibitem{Sezgin}  A. Sezgin, A.O. Atag\"{u}n,      \emph{On operations of soft sets},
Computers and Mathematics with Applications, 61(2011), 1457-1467.

\bibitem{Singh}   A.K. Singh, P. Parida,      \emph{Soft Computing in Financial Decision Making},
Global Journal of Management and Business Studies, 3/2(2013), 103-110.

\bibitem{Vijayabalaji}   S. Vijayabalaji,  A. Ramesh,      \emph{A New Decision Making Theory in Soft Matrices},
International Journal of Pure and Applied Mathematics, 86/3(2013), 927-939.

\bibitem{Zhang}   Z. Zhang, C.   Wang, D. Tian, K. Li,       \emph{A novel approach to interval-valued intuitionistic fuzzy soft set based decision making},
Applied Mathematical Modelling, In Press(2013).

\end{thebibliography}
\end{document}